\documentclass[12pt,a4paper]{amsart}
\usepackage{amssymb,mathrsfs, graphicx, subfig, wasysym}
\usepackage{hyperref}
\usepackage[german, english]{babel}
\newtheorem{theorem}{Theorem}
\newtheorem{lemma}{Lemma}[section]

\theoremstyle{definition}

\newtheorem{example}[lemma]{Example}

\theoremstyle{remark}

\numberwithin{equation}{section}

\newcommand{\N}{{\mathbb N}}
\newcommand{\T}{{\mathbb T}}

\newcommand{\A}{{\mathcal A}}

\newcommand{\CC}{{\mathbb C}}
\newcommand{\NN}{{\mathbb N}}
\newcommand{\TT}{{\mathbb T}}
\newcommand{\PP}{{\mathbb P}}
\newcommand{\EE}{{\mathbb E}}

\newcommand{\HH}{{\mathbb H}}
\newcommand{\h}{{\mathbb H}}

\newcommand{\de}{{\delta}}

\newcommand{\ph}{{\varphi}}
\newcommand{\lm}{{\lambda}}
\newcommand{\gm}{{\gamma}}
\newcommand{\si}{{\sigma}}

\newcommand{\Gm}{{\Gamma}}

\newcommand{\Hm}[1]{\leavevmode{\marginpar{\tiny%
$\hbox to 0mm{\hspace*{-0.5mm}$\leftarrow$\hss}%
\vcenter{\vrule depth 0.1mm height 0.1mm width \the\marginparwidth}%
\hbox to 0mm{\hss$\rightarrow$\hspace*{-0.5mm}}$\\\relax\raggedright
#1}}}

\newtheorem*{coron}{\textbf{\rm Corollary}}

\begin{document}
\title[Trees of finite cone type]{An invitation to trees of finite cone type: random and deterministic operators}

\author{Matthias Keller}
\address{Mathematisches Institut \\
Friedrich Schiller Universit{\"a}t Jena \\
07743 Jena, Germany } \email{m.keller@uni-jena.de}

\author{Daniel Lenz}
\address{Mathematisches Institut \\
Friedrich Schiller Universit{\"a}t Jena \\
07743 Jena, Germany } \email{daniel.lenz@uni-jena.de}

\author{Simone Warzel }
\address{Technische Universit\"at M\"unchen
Zentrum Mathematik\\  Boltzmannstraße 3, 85747 Garching   }
\email{warzel@ma.tum.de}

\subjclass[2000]{}
\date{\today}

\parindent0cm

\begin{abstract}
Trees of finite cone type have appeared in various contexts. In
particular, they come up as simplified models of regular
tessellations of the hyperbolic plane. The spectral theory of the
associated Laplacians can thus be seen as induced by geometry. Here
we give an introduction focusing on background and then turn to
recent results for (random) perturbations of trees of finite cone
type and their spectral theory.
\end{abstract}

\maketitle

{\small\tableofcontents}

\section{Introduction}
Trees of finite cone type are a generalization of regular trees.
They have a special form of 'recursive structure'. Over the years
they  have appeared in various contexts. These include -- in the
order of appearance -- the following:

\begin{itemize}

\item Theory of random walks.

\item Spectral theory of hyperbolic tessellations.

\item Random Schr\"odinger operators.

\end{itemize}

Here, we want to give an introduction into the topic. In particular,
we want to illuminate each of these contexts and provide some
background. Our main results then concern (random) Schr\"odinger
operators associated to such trees.

More specifically, this article will deal with the following
questions in the indicated sections.

\begin{itemize}
\item What are trees of finite cone type? (Section \ref{What}).
\item Why study trees of finite cone type? (Section \ref{Why}).
\item Which spectral properties have trees of finite cone type?
(Section~\ref{Which}).
\item How are trees of finite cone type related to multi-type Galton Watson trees? (Section~\ref{How}).
\end{itemize}
Finally, we will present some ideas behind the proofs in Section
\ref{Some}.

\smallskip

The results in Section~\ref{Which} are published in \cite{KLW,KLW2}.
These works in turn rely on   the PhD thesis of one of the authors
(M.K.). The results of Section~\ref{How} for Galton Watson type trees can be found in
\cite{Kel2}.

\medskip

\textbf{Acknowledgments.} D.L. and S.W. take this opportunity to
express their heartfelt thanks to  the organizers of Pasturfest
(2013) in Hagen. Furthermore, M.K. would like to thank Balint Virag
for the very inspiring discussions during his visit in Toronto
leading to simplifications of the proofs presented in this survey.

\section{What are trees of finite cone type?}\label{What}
In this section we introduce trees of finite cone type.  We will
assume familiarity with basic notions from graph theory and just
remind the reader of certain notions.

\medskip

For us a graph is a pair $(V,E)$ consisting of a set $V$ called
\textit{vertices} and a set $E$ called \textit{edges} consisting of
subsets of $V$ with exactly two elements. For $x,y\in V$,  we will
write $x\sim y$ if $\{x,y\} \in E$ and call $x$ and $y$ neighbors.
The number of neighbors of $x$ is called the \textit{degree} of~$x$.

A finite sequence $x_0,x_1,\ldots, x_n$ of pairwise different
vertices is called a \textit{path of length $n$ connecting $x_0$ and
$x_n$} if $x_j \sim x_{j+1}$ holds  for all $j=0,\ldots, n -1$. If
any two vertices are \emph{connected} by a path, then the graph is
called connected. In this case, the distance between two vertices is
the length of the shortest path connecting these vertices. If for
any two vertices there exists exactly one path connecting them (thus
the graph is connected without loops), then the graph is called a
\textit{tree} and the unique path between two vertices is called a
\textit{geodesic}. A connected graph with a distinguished vertex is
called a \textit{rooted graph} and the distinguished vertex is
called the \textit{root}.

In a rooted tree (i.e., a tree with a root) the elements with
distance $n$ to the root are called $n$-sphere. A neighbor of a
vertex $v$ in the $n$-sphere is called a \textit{forward neighbor}
if it belongs to the $(n+1)$-sphere. The \textit{cone}, or rather
\textit{forward cone}, of a vertex $v$ in a rooted tree is then the
smallest subgraph of the rooted tree which contains $v$ and all
forward neighbors of any of its vertices. A rooted tree is called
\emph{regular} if the number of forward neighbors of any vertex is
constant.

We will deal with graphs in which every vertex is labeled. This
means we will have a map from the vertices to some (finite) set
called the labels.

\medskip

Trees of finite cone type arise from the following two  pieces of
data:

\begin{itemize}
\item[(D1)] A finite set $\mathcal{A}$  called \textit{vertex labels}.

\item[(D2)] A map $M : \mathcal{A} \times \mathcal{A}  \longrightarrow
\NN_0$ called \textit{substitution rule}.

\end{itemize}

\smallskip

Given these  data we can then construct a rooted   tree $\TT =
\TT(M,j)$ whose vertices carry labels from $\mathcal{A}$ and whose
root is labeled by $j$ in the following way:

\medskip

\textit{Every vertex with label $k$ of the $n$-sphere is joined to
$M_{l,k}$ vertices of label $l$ of the $(n+1)$-sphere.}

\medskip

The arising structure is obviously a rooted tree  in which every
vertex has a label. By construction the cone of any vertex is
completely determined by the label of the vertex. In particular,
there are only finitely many different types of cones present.
Conversely, it is not hard to see that any rooted tree with the
property that it has only finitely many types of cones arises in the
manner described above. We will refer to such trees as \textit{trees
of finite cone type}.  By construction the vertex degree  is
uniformly bounded in a tree of finite cone type.

Note that our data  allows us to  constructed one tree for each
element of $\mathcal{A}$.

\medskip

\textbf{Remark.} As already mentioned in the introduction, trees of
finite cone type have come up in several contexts. In fact, they
have also been denoted as \textit{periodic trees} or as
\textit{monoids} with certain properties. Here, we follow  Nagnibeda
/ Woess \cite{NW} in calling them trees of finite cone type.

\medskip

In order to achieve a meaningful theory we will impose the following
\textbf{assumptions} on our data or more precisely on the
substitution matrix $M$:

\begin{itemize}
\item[(M0)] If $\mathcal{A}$ consists of only one point, then $M = M_{1,1}\geq 2$ (\textit{not one-dimensional}),
\item[(M1)] $M_{j,j}\ge1$ for all $j\in \A$ (\textit{positive diagonal}),
\item[(M2)] there exists $n = n (M)$ with $M^n$ has positive entries  (\textit{primitivity}).
\end{itemize}

\medskip

\textbf{Remark.} Due to positivity of diagonal primitivity follows
already from a generally weaker condition known as  irreducibility.

\medskip

Let us discuss a concrete example, which we call the
\textit{Fibonacci tree}.
\begin{example}
Let $\A=\{\circ,\bullet\}$ and $M=\left(
\begin{array}{cc}
2 & 1 \\
1 & 1 \\
\end{array}
\right) $. The tree $\T=\T(M,2)$ is illustrated in
Figure~\ref{f:tree} from \cite{KLW}. Obviously, the sequence of Fibonacci numbers
arises if  we count the number of vertices with label $\bullet$ and
with label  $\circ$  respectively in each sphere.

\begin{figure}[!h]
\centering
\scalebox{0.3}{\includegraphics{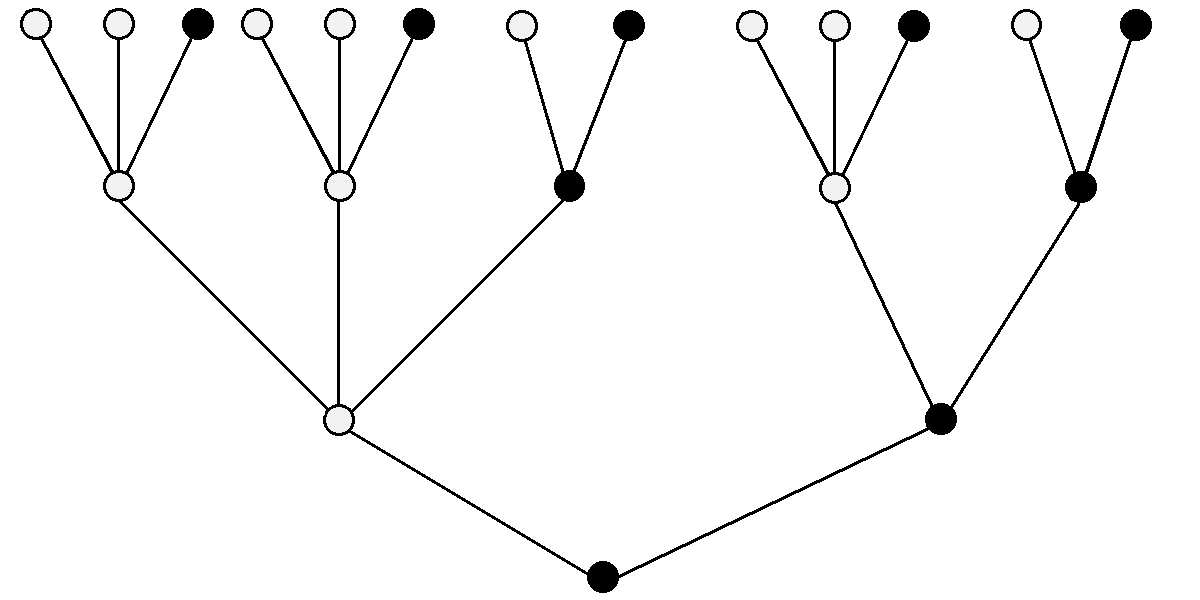}}
\caption{An example of a tree constructed by a substitution matrix.}\label{f:tree}
\end{figure}
\end{example}

\medskip

After we have introduced the graphs in question,  we will now come to
the associated operators. Given a tree $\TT = \TT (M,j)$ as above
with vertex set $V$ we consider the Hilbert space
$$\ell^2 (\TT):=\ell^2 (V) :=\{ f : V\longrightarrow \CC :
\sum_{v\in V} |f(v)|^2 < \infty\}$$ with inner product given by
$$\langle f, g\rangle = \sum_{v\in V} \overline{f (v) } g(v).$$

In our context the analogue  of  the usual Laplacian  in the
Euclidean situation  will be given by the \textit{adjacency matrix}
$$ L: \ell^2 (\TT) \longrightarrow \ell^2 (\TT)$$
$$ (L \varphi) (x) = \sum_{y\sim x} \varphi (y).$$
The  operator $L$ can easily be seen to be bounded (as the vertex
degree is uniformly bounded). Moreover, the operator $L$ is clearly
symmetric. Hence, $L$ is a selfadjoint operator. We will be
interested in the spectral theory of $L$ and later of (random)
perturbations of it.

\medskip

\textbf{Remark.} Various extension of the operator $L$ can be
treated by our methods. We could, for example,  include edge weights
and / or potentials depending only on the labeling. Also, it is
possible to replace  the underlying (constant)  measure on the set
of vertices by any label invariant measure. Specifically, all our
results below remain valid for the operator
$$ \Delta  := \mbox{I } - D^{-1} L$$
on the weighted $\ell^2$ space $\ell^2 (V, D)$. Here,
 $I$ is the identity and $D$ is the operator of multiplication with the
degree. This operator is sometimes known as normalized Laplacian. It
is a selfadjoint operator and  occurs in the investigation of random
walks (see below).  We refrain from further discussing this here in
order to simplify the presentation, see \cite{Kel2} for details.

\section{Why study  trees of finite cone type?}\label{Why}
In this section we  provide some background on trees of finite cone
type. In particular, we discuss the following three contexts in
which these graphs and their associated operators came up.

\begin{enumerate}
\item Random operators on trees.
\item Simplified models of regular tessellations of the hyperbolic plane.
\item Random walks on trees.
\end{enumerate}

\subsection{Random operators on trees}
After the seminal work of the theoretical physicist Anderson on
random media in the 50ies, the mathematically  rigorous
investigation of random operators gained momentum  in the 70ies
beginning with work  due to mainly  Pastur,  Molchanov, Goldshtein.
A main focus of the investigations is the phenomenon of
localization, i.e., intervals with pure point spectrum with
exponentially decaying eigenfunctions. In the one-dimensional
situation a most general result of  Carmona / Klein / Martinelli
\cite{CKM} shows localization on the whole spectrum. For the
treatment of models in arbitrary dimensions key results were
provided by the works of Fr\"ohlich / Spencer \cite{FS} and
Fr\"ohlich / Martinelli / Scoppola / Spencer, and von Dreyfuss /
Klein \cite{vDK}. This lead to an approach called multi-scale
analysis. This was later complemented by a different approach known
as Aizenman-Molchanov method \cite{AM}. Building on these approaches
a large body of work has accumulated over the last three decades and
it  seems fair to say that localization is quite well understood by
now.

The big open question now concerns 'existence of  extended states',
i.e., the occurrence of an absolutely continuous component in the
spectrum. For one dimensional models such a component is known to be
absent. In fact, the  spectrum is known to be pure point (see
above). In two dimensional models the  behavior expected by
physicists is under dispute. But for higher dimensions, i.e., from
dimension 3 on, the occurrence of an absolutely continuous component
in the spectrum of a random model is conjectured to take place.  So
far, this could not be proven.

On the other hand, work of Klein \cite{Kle,Kle2} in the 90ies showed
that  an absolutely continuous component exists in a perturbative regime for infinite
dimensional random models, i.e., for models on trees. So, random models
on trees are the only models among the 'usual' random operators where the extended states conjecture is proved. Various groups have taken up this
line of investigation. In fact,
both \cite{ASW}  and \cite{FHS,FHS2} come up with alternative approaches at about
the same time.

Beyond the perturbative regime, some surprises have been discovered concerning the regime of extended states on regular trees. Among them is that even at weak disorder, the extended states can be found well beyond energies of the unperturbed model into the regime of Lifshitz tails. The mechanism for the appearance of extended states in this non-perturbative regime are disorder-induced resonances \cite{AW}.

\begin{figure}[!h]
\centering
\scalebox{0.3}{\includegraphics{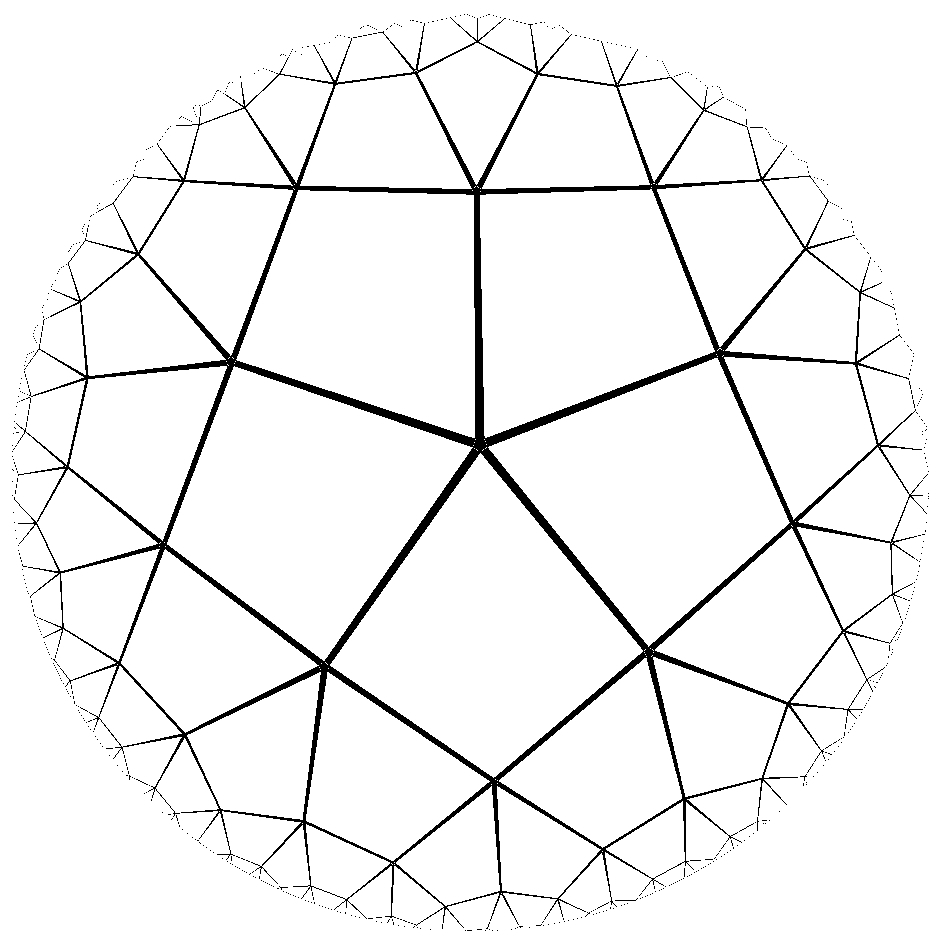}}
\scalebox{0.3}{\includegraphics{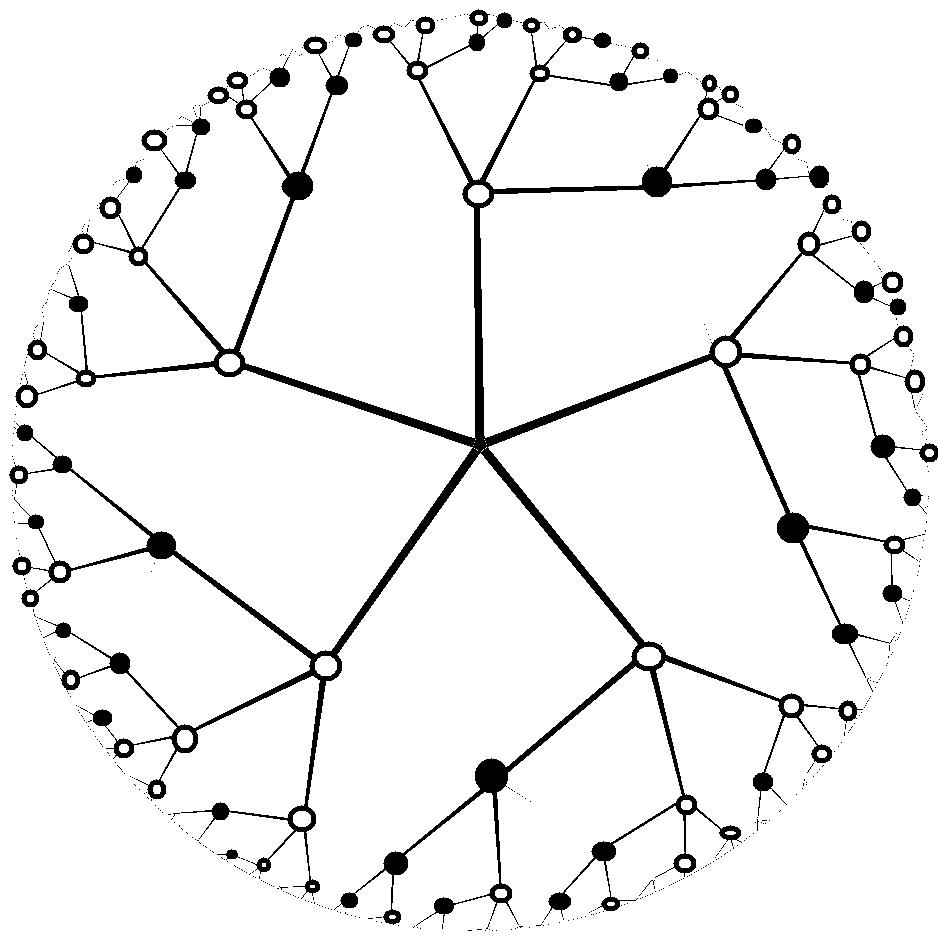}}
\caption{The regular hyperbolic (4,5)-tessellation and a spanning tree which is constructed by a substitution matrix.}
\label{f:tess}
\end{figure}

All these works deal with regular trees.
This restriction to regular trees is rather by convenience than by
intrinsic (physical) reasons. Thus, it is reasonable to consider more
complicated models. In fact,  the authors of the present article
started their corresponding  investigations with an attempt to study
spectral theory of regular tessellations of the hyperbolic plane. The
tessellation models turned out to be too hard to analyze due to the
appearance of loops in the graphs. Then, systematic ways of cutting
loops lead us to models of trees of finite cone type (compare
\cite{Kel2} for details), see Figure~\ref{f:tess} from \cite{Kel1}.

A  similar occurrence of trees of finite
cone type in the spectral theory of regular tessellations of the
hyperbolic plane will be discussed below.

As for recent works on various variants of trees we also mention
\cite{FHS3,FHH,Hal,KS,Sha}

\subsection{Regular tilings of the hyperbolic plane}
Trees of finite cone type have also come up in the study of spectral
theory of tilings of the hyperbolic plane. This will be discussed in
this section.

\medskip

The three basic models of planar geometry are the Euclidean plane
(with curvature zero), the sphere (with curvature 1) and the
hyperbolic plane (with curvature $-1$). Each of these three models
allows for regular tessellations and the spectral theory of the
(discrete) Laplacians associated to these tessellations has been
quite a focus of research. Concerning this spectral theory of
regular tessellations the basic picture is as follows:

\begin{itemize}
\item Euclidean plane: band spectrum, purely absolutely continuous
spectrum.
\item Sphere in three dimensional  Euclidean space:  finite pure point spectrum
\item Hyperbolic plane: quite unknown!
\end{itemize}
Let us comment on the last point. Spectral theory may be divided
into investigation of the form of the spectrum and investigation of
the type of the spectrum. For both parts only  some pieces of
information are known for regular tessellations of the hyperbolic
plane.

Quite some effort has been put into studying the infimum of the
spectrum.\footnote{Our subsequent discussion does not refer to the
operator $L$ but rather to its normalized  version $\Delta =I -
D^{-1} L$ on $\ell^2 (V,D)$ (see above).}
 It is known that this infimum is positive
(due to basic Cheeger estimates). The exact value is unknown. In
fact, there is a famous conjecture of Sarnak stating that the
infimum of the spectrum is an algebraic number (see the book
\cite{Woe} of Woess for discussion). While the conjecture is open,
some  works have been devoted to giving  estimates on the actual
value of the infimum of the spectrum. Among these we mention
\cite{BCCSH,BCS,CV,Nag,Zuk}.  In particular, the  best lower bounds
known so far  is  established  by Nagnibeda  in \cite{Nag} and the
best upper bound by Bartholdi / Ceccherini-Silberstein in
\cite{BCS}. Quite remarkably these bounds agree up to the third
digit. In a very recent preprint, the upper bounds of \cite{BCS}
seem to have been improved once more by Gouezel in \cite{Gou}. (Note
that these works use a somewhat different normalization and thus
their use of lower and upper bound is just the reverse of our use
here.)

In  our context the mentioned work of Nagnibeda is most relevant. In
fact, her approach to the lower bound can be divided in three steps:
In the first step she associates a tree of finite cone type to a
given tessellation. This tree arises as some sort of universal
covering of the tessellation. In the second step she shows that the
infimum of the spectrum of the tessellation can be bounded below by
the infimum of the spectrum of the tree. In the last  step she then
bounds the spectrum of the tree. So, trees of finite cone type most
naturally come up in the study of spectral theory of regular
tessellations.

As for further results on spectral theory it seems that for those
regular tessellations of the hyperbolic plane which are Cayley graphs
of torsion free groups some information on both parts of spectral
theory is available. More specifically, after talking with various
experts  it seems to us that for such models there should be  no
eigenvalues as a consequence of  (an established version of) Linnells
zero-divisor conjecture and that the spectrum should be an interval as
a consequence of (an established version of) Baum-Connes conjecture.
It is unclear (to us at least) whether explicit statements to these
effects can be found in the literature. In any case, no information
on existence of absolutely continuous spectrum for such models seems
available at the moment.

For those regular tessellations of the hyperbolic plane which are not
Cayley graphs of torsion free subgroups even less is known. It is
tempting to think of these models as certain forms of coverings of a
torsion free situation. This suggests that they  have bands of absolutely continuous spectrum and possibly finitely many  eigenvalues in between the bands. The extend to which such a
reasoning could be made precise is unclear at the moment. Two of the
authors of the present paper (M.K. and D.L.)  are currently working
together with Lukasz Grabowski on this topic, \cite{GKL}.

\medskip

\textbf{Remark.} To us it seems most remarkable that so little is
known on spectral types of (discrete) Laplacians on  regular
tessellation of the hyperbolic plane. This is somewhat reminiscent
of the situation encountered in the study of quasicrystals. There
essentially nothing is known on the spectral type of such basic
models as Laplacian associated to Penrose tiling. In both cases once
encounters very natural geometric situations exhibiting a large
degree of overall order. Still it seems hard to get a grip on. On
the technical level this seems to be related to the lack of suitable
Fourier type analysis.

\subsection{Random walk on trees of finite cone type}
Trees of finite cone type have also been studied in the context of
random walks. This will be briefly discussed in this section.

\medskip

Let $\TT$ be  tree of finite cone type and denote the operator of
multiplication by the degree by $D$. Then, $P=D^{-1}L = I - \Delta $
is a selfadjoint operator on the Hilbert space with the vertex
degree as a measure. It can be seen as the transition matrix of a
random walk (in discrete time). Properties of this random walk have
been studied in particular in \cite{Lyo,Tak,NW,Mai}.  These works
are mostly concerned with recurrence and transience, i.e.,  -- very
roughly speaking -- the behavior of the Greens function at the
infimum of the spectrum.

\section{Which  spectral properties have  trees of finite cone type?} \label{Which} In this section we present some results on
spectral theory of trees of finite cone type.

\medskip

To be specific we will consider the following {situation:}
$\mathcal{A}$  is the set of vertex labels,  $M : \mathcal{A} \times
\mathcal{A}  \longrightarrow \NN_0$ is the substitution matrix and
the tree $\TT = \TT(M,j)$  with root labeled by $j$ is generated
according to the rule given above that every vertex with label $k$
of the $n$-th sphere is joined to $M_{l,k}$ vertices of label $l$ of
the $(n+1)$-sphere. Moreover, we will assume that the assumptions
(M0), (M1) and (M2) are in effect. We will then consider the
operator $ L: \ell^2 (\TT) \longrightarrow \ell^2 (\TT)$, $ (L
\varphi) (x) = \sum_{y\sim x} \varphi (y)$ introduced above.

\subsection{The unperturbed operator}
In this section we present the basic result on the spectral theory
of the unperturbed operator.

\bigskip

\begin{theorem}[\cite{KLW}] There exist finitely many intervals such that for every $j\in \mathcal{A}$ the spectrum of $L$ associated to $\TT(M,j)$ consists of exactly these intervals and is purely absolutely continuous.
\end{theorem}

\textbf{Remarks.} (a) There is an  effective bound on the number of
bands by results of Milnor \cite{Mil}  on algebraic varieties.

(b)  A tree constructed from $M$ not satisfying $M_{j,j}>0$ may have
eigenvalues.

(c) The result  generalizes the well-known  spectral theory of
regular trees.

\subsection{Radially symmetric perturbations}

A function $v$ on the vertices is  \emph{radially labeling symmetric} if  vertices with the same label in the same generation yield the same value.

\bigskip

\begin{theorem}[\cite{KLW}]  Assume $\TT$ is not a regular tree. Then, for every compact $I$ contained in the interior of  $\sigma (L) \setminus \{0\}$, there exists $\lambda>0$ such that  for all radially labeling symmetric $v $ with values in $[-1,1]$ we have
$$I \subset \sigma_{ac} ( L  + \lambda v)\:\; \mbox{and}\:\; I \cap \sigma_{sing} ( L + \lambda v) = \emptyset.$$
\end{theorem}

\textbf{Remarks.}   The  corresponding statement  is  quite invalid
for regular trees. There one is essentially in a one-dimensional
situation (see e.g. \cite{Bre}) and then a random potential yields
pure point spectrum.  This is already discussed in \cite{ASW}  (see
\cite{BF} as well for related material).

\bigskip

\begin{coron}[\cite{KLW}]\textit{  Assume $\TT$ is not a regular tree. Then, for every radially labeling symmetric $v$ with $v(x) \to 0$ for $x\to \infty$, we have
$$\sigma_{ac} (L) =  \sigma_{ac} ( L  + v).$$}
\end{coron}

\textbf{Remark.} For regular trees  one is far from having such a
general statement compare e.g. \cite{Kup}.

\subsection{Random perturbations}
Here, we consider random perturbations of our model. Thus, beyond
the operator $L$ as above we are given the following pieces of data:

A probability space $(\varOmega,\PP)$ and a measurable  map
$$ (v,\vartheta) : \varOmega \times V\longrightarrow (-1,1)\times (-1,1),\:\; (\omega, x)\mapsto (v_x^\omega, \vartheta_x^\omega)$$
satisfying the following two properties:

\begin{itemize}

\item For all $x,y\in V$ the random variables  $(v_x,\vartheta_x)$ and $(v_y,\vartheta_y)$ are independent if the forward trees of $x$ and $y$ do not intersect.

\item For all $x,y\in V$ that share the same label the restrictions of the random variables $(v,\vartheta)$ to the isomorphic forward trees of $x$ and $y$ are identically distributed.
\end{itemize}
For a vertex $x$ that is not the root let $\dot x$ be the vertex such that $x$ is a forward neighbor of $\dot x$ and denote
\begin{align*}
    \vartheta_{x,\dot x}=\vartheta_{\dot x,x}=\vartheta_{x}.
\end{align*}
We then consider the operator

$$(H^{\lambda,\omega} \varphi) (x) = \sum_{y\sim x} (\varphi (y) + \lambda \vartheta^{\omega}_{x,y}\varphi (y)) + \lambda v_x^\omega \varphi (x).$$

\medskip

\bigskip

\begin{theorem}
There exists a finite set $\Sigma_0 \subset \sigma (L)$ such that
for each compact $I\subset \sigma (L)\setminus \Sigma_0$ there is
$\lambda>0$ such that  $H^{\lambda,\omega} $ has almost surely
purely absolutely continuous spectrum in $I$ for all $(v,\vartheta)$
satisfying the conditions above.
\end{theorem}

\textbf{Remarks.} (a) The theorem includes regular trees as a
special case.

(b)   If the tree is not regular (and in some other cases as well),
$\Sigma_0$ is just the set of boundary bounds of $\sigma (L)$
together with $0$.

\section{How are trees of finite cone type related to multi-type Galton Watson trees?}\label{How}

In this section we relate trees of finite cone type to multi-type Galton Watson trees. We present a result stating that whenever certain Galton Watson trees are close to a tree of finite cone type in distribution, then these random trees inherit some of the absolutely continuous spectrum from the tree of finite cone type. The operators on such trees are also random operators, where now the randomness is induced by the underlying combinatorial geometry rather than by a random perturbation of the parameters.

A multi-type Galton Watson branching process $b$ with finite set $\A$ of \emph{types} or \emph{labels} is determined by  numbers $P_{ij}^{(b)}(n)$  such that $\sum_{n=0}^{\infty}P_{ij}^{(b)}(n)=1$, $i,j\in\A$. The number $P_{ij}(n)$ encodes the probability that a vertex with label $i$ has $n$ forward neighbors of label $j$.
We may encode the possible configurations  of forward neighbors of a vertex by vectors $s\in \N_{0}^{\A}$ (i.e., the vertex has $s_{j}$ forward neighbors of type $j$) and define by
\begin{align*}
    \PP_{i}^{(b)}(s)=\sum_{j\in \A}P_{ij}^{(b)}(n)
\end{align*}
the probability that a vertex of type $i\in\A$ has the forward configuration $s$. Moreover,  $\|s\|_{1}=\sum_{j\in A}s_{j}$ denotes the total number of forward neighbors of a configuration $s$.

We impose two further assumptions on the branching processes. The first is a mild growth condition and the second says that every vertex has at least one forward neighbor:
\begin{itemize}
  \item [(B1)] There is $p>1$ such that  $\sum_{s\in \N_{0}^{\A}}\PP^{(b)}_{j}(s)\|s\|_{1}^{p}<\infty$ for all $j\in \A$.
  \item [(B2)] $\PP^{(b)}_{j}(0)=0$ for all $j\in \A$.
\end{itemize}

We observe that a tree of finite cone type given by a labeling set $\A$ and a substitution matrix $M$ can be viewed as such a multi-type Galton Watson branching process $b_{M}$ by letting $P_{ij}^{b_{M}}(n)=1$ for $n=M_{i,j}$ and $0$ otherwise. Hence,
\begin{align*}
    \PP^{(b_{M})}_{i}(M_{i,\cdot})=1\quad\mbox{and}\quad     \PP^{(b_{M})}_{i}(s)=0 \mbox{ for $s\neq M_{i,\cdot}$}.
\end{align*}
Obviously, (B1) is satisfied and (B2) follows from (M1).  In this sense, trees of finite cone type are multi-type Galton Watson trees with deterministic distribution.

We introduce a metric $d_{p}$ on branching processes with the same type set $\A$ that satisfy (B1) for the same $p>1$ via
\begin{align*}
    d_{p}(b_{1},b_{2})=\max_{i\in \A}\sum_{s\in \N_{0}^{\A}} |\PP_{i}^{(b_{1})}(s)-\PP_{i}^{(b_{2})}(s)|\|s\|_{1}^{p}.
\end{align*}

We denote the set of realizations of a Galton Watson branching process by $\Theta(b)$. For $\theta=(V^{\theta},E^{\theta})\in \Theta(b)$  we denote the adjacency matrix on $\ell^{2}(\theta)=\ell^{2}(V^{\theta})$ by $L^{\theta}$. The operator $L^{\theta}$ may not be bounded, but by techniques of \cite{GS} its restrictions to the functions of finite support is almost sure essentially selfadjoint on $\ell^{2}(\theta)$.

Given a fixed tree of finite cone type $\TT$ we can now state a theorem proven in \cite{Kel2}. Loosely speaking it says that operators $L^{\theta}$ on random trees $\theta$ close to $\TT$ in distribution inherit almost surely some of the absolutely continuous spectrum of $L$.

\begin{theorem} There exists a finite set $\Sigma_0 \subset \sigma (L)$ such that for each compact $I\subset \sigma (L)\setminus \Sigma_0$ there is $\lambda$ such that if a multi-type Galton Watson branching process $b$ satisfies (B1), (B2) and
\begin{align*}
     d_{p}(b,b_{M})<\lm,
\end{align*}
then
$L^{\theta}$ has  purely absolutely continuous spectrum in $I$ for almost all $\theta\in \Theta(b)$.
\end{theorem}
A strategy for a proof of such a theorem for certain single-type Galton Watson trees was also outlined in \cite{FHS3}.

\section{Some ideas behind the proofs}\label{Some}

For a tree $\TT=\TT(M,j)$, $j\in \A$, and a vertex $o$ let $\TT_{o}$ be the forward tree with root $o$.
The  forward Green function  of an operator $H$, acting as
\begin{align*}
    H\ph(x)=\sum_{y\sim x}t_{xy}\ph(y)+v(x)\ph(x)
\end{align*}
on $\ell^{2}(\TT)$,  at $o$ is defined by
$$\varGamma_{o} (z,H):= \langle ( H_{\TT_{o}} - z)^{-1} \delta_o, \delta_o\rangle ,\quad z\in \HH,$$
where $\HH=\{z\in \CC\mid \Im z>0\}$. The second resolvent formula gives the \emph{recursion formula}
$$ - \frac{1}{\varGamma_{o} (z,H)} =  z  -v(o) + \sum_{y\in S_{o}} |t_{oy}|^{2} \varGamma_y (z,H),$$
where $S_{o}$ is the forward sphere of $o$.

\subsection{The unperturbed operator}
We present a basic inequality which shows the absolute continuity of the spectral measure  of the unperturbed operator $L$ at the root.

For $H=L$ the forward Green function $\Gm_{x}(z,L)$ is equal for all vertices $x$ with  the same label $j$, thus,  we write
$\Gm_{j}(z,L)$. The recursion formula now reads
$$ - \frac{1}{\varGamma_{j} (z,L)} =  z   + \sum_{k\in\A} M_{j,k} \varGamma_k (z,L),$$

Taking imaginary parts in the recursion formula and multiplying by $|\Gm_{j} (z,L)|^{2}$,
we arrive at
\begin{align*}
\Im\Gm_{j} (z,L)=\sum_{k\in\A}{M_{j,k}}\Im \Gm_{k} (z,L)|\Gm_{j} (z,L)|^{2} \ge M_{j,j}\Im\Gm_{j} (z,L)|\Gm_{j} (z,L)|^{2}
\end{align*}
since $\Gm_{k} (z,L)$ are Herglotz functions and, therefore, $\Im\Gm_{k} (z,L)>0$ for $z\in\h$.
Hence,
\begin{align*}
    |\Gm_{j} (z,L)|\leq{{|M_{j,j}|}^{-\frac{1}{2}}}.
\end{align*}
This gives that the Green function is uniformly bounded on $\HH$. By
standard results,  the spectral measure associated to the  delta
function of the root must then be  absolutely continuous. The result
for arbitrary spectral measures follows by further considerations
using the recursion formula, confer \cite[Proposition~3]{KLW}.

\subsection{Random operators}
In this section we discuss the strategy of the proof of absolutely continuous spectrum for random operators.
Consider the function $\gm:\HH\times\HH\to[0,\infty)$
\begin{align*}
    \gm(g,h)=\frac{|g-h|^{2}}{\Im g\Im h},\qquad g,h\in \HH,
\end{align*}
which is related  to the hyperbolic metric $d_{\HH}$ of the upper half plane viz
\begin{align*}
    d_{\HH}=\cosh^{-1}(\tfrac{1}{2}\gm+1).
\end{align*}
The function $\gm$ is symmetric and zero only at the diagonal, but it does not satisfy the triangle inequality. Nevertheless, $\gm$ works very well with the recursion formula.

The recursion formula  induces a map $\Phi:\h^{S_{x}}\to\h$
\begin{align*}
(g_{x})_{x\in S_{o}}\mapsto -\frac{1}{z-v(o)+\sum_{x\in S_{o}}g_{x}},
\end{align*}
mapping $(\Gm_{x}(z,H))_{x\in S_{o}}$ to $\Gm_{o}(z,H)$,
which decomposes into maps $$\rho\circ\si\circ\tau:\h^{S_{o}}\longrightarrow\h\longrightarrow\h\longrightarrow\h$$
given by
\begin{align*}
    \tau((g_{x}))= \sum_{x\in S_{o}}g_{x},\quad \si(g)= z-v(o)+g,\quad \rho(g)=-\frac{1}{g}.
\end{align*}

The map $\Phi$ turns out to be a hyperbolic contraction in a certain sense, e.g. $\gm(\Phi(g),\Phi(h))<\sup_{y\in S_{x}}\gm(g_{y},h_{y})$, $g,h\in\h^{S_{o}}$. Indeed,
$\rho$ is an isometry and $\si$ is a uniform contraction, however, 'unfortunately' with the contraction coefficient going to one as on $\Im z\to0$. Finally, the map $\tau$ can be seen to be a hyperbolic mean, (i.e., as $g\mapsto |S_{x}|g$ is a hyperbolic isometry we may also consider $\widetilde\tau((g_{x}))=\frac{1}{|S_{x}|}\sum_{y}g_{y}$ instead of $\tau$). A hyperbolic mean may contract the input values for  two reasons: either two of the vectors $g_{x}-h_{x}$, $x\in S_{o}$, are  linearly independent or two of the quantities $|g_{x}-h_{x}|/\Im g_{x}\Im h_{x}$, $x\in S_{o}$ are different. This is exploited in the proof.

We further use  that the random variables $\Gm_{x}(z,H^{\lm,\omega})$,  are independent for $x\in S_{o}$, and identically distributed for all $x$ with the same label. Consequently, we introduce the random variables $G:\Omega\to[0,\infty)^{\A}$
\begin{align*}
    G_{j}=\gm (\Gm_{x}(z,H^{\lm,\omega}),\Gm_{x}(z,L))^{p}
\end{align*}
for some $p>1$ and $x $ having the label $j$.

The crucial estimate is then the following vector inequality.
We show that there are $\de>0$, a stochastic $(|\A|\times|\A|)$ matrix $P$ and $C$ such that
\begin{align*}
    \EE(G)\leq (1-\de)P    \EE(G)+C,
\end{align*}
where the inequality is of course interpreted componentwise.

By the Perron-Frobenius theorem there is a positive left eigenvector $u$ of $P$ and we arrive at
\begin{align*}
\langle u,G\rangle\leq(1-\de) \langle u,PG\rangle+ C=(1-\de)\langle u,G\rangle+C
\end{align*}
and, hence,
\begin{align*}
\EE(    \gm (\Gm_{x}(z,H^{\lm,\omega}),\Gm_{x}(z,L))^{p})=\EE(G_{j})\leq \frac{C}{\de u_{j}}
\end{align*}
for all $x$ where $j$ is the label of $x$.

By standard arguments found in  \cite{Kle2}, \cite{FHS2} and \cite{Kel1,KLW2}, one deduces
\begin{align*}
    \liminf_{\eta\to0}\int_{I} |\langle(H^{\lm,\omega}-E-i\eta)^{-1}\de_{x},\de_{x}\rangle |^{p} dE<\infty
\end{align*}
almost surely. This, however, yields purely absolutely continuous spectrum in $I$ by a limiting absorption principle, see \cite{Kle2}.

At the end we indicate how the proof for Galton Watson trees fits in the scheme above. For realizations $\theta\in\Theta(b)$. The probability that the first two spheres in $\theta$ look exactly like the ones of the tree of finite cone type is very large by the assumption $d_{p}(b,b_{M})<\lm$.

To prove the crucial vector inequality above, we divide the expected value in the 'good part', where the spheres are exactly the same as in the deterministic case, and the 'other part'. It turns out that the error made in the 'other part' is bounded which is then swallowed by the small probability of its occurrence. For the 'good part' we  apply the proof for random potentials above to conclude the statement.

\end{document}